\documentstyle[12pt]{amsart}
\theoremstyle{plain}
\newtheorem{Thm}{}[section]
\title{ Injectivity of the symmetric map for line bundles}
\author{montserrat teixidor i bigas}
\address {Mathematics Department, Tufts University, Medford,
MA02155}
 \email{montserrat.teixidoribigas@@tufts.edu}
\begin{document}
\maketitle

\begin{abstract}
Let $C$ be a generic non-singular curve of genus $g$ defined over
a field of characteristic different from two. We show that for
every line bundle on $C$ of degree at most $g+1$, the natural
product map
$$S^2(H^0(C,L))\rightarrow H^0(C,L^2)$$ is injective. We also show
that the bound on the degree of $L$ is sharp
\end{abstract}

\begin{section}{Introduction}

Linear series have played a key role in the study of curves (and
higher dimensional varieties). It is usually important to know how
the spaces of sections of line bundles relate to each other. For
example for a generic curve of genus $g$, the Petri map
$$H^0(C,L)\otimes H^0(C,K\otimes L^{-1})\rightarrow H^0(C,K)$$
is injective. This property gives a lot of information on the
structure of the set of linear series of a given dimension. On the
other hand, the surjectivity of the maps
$$S^n(H^0(C,L))\rightarrow H^0(C,L^n)$$
ensure the projective normality of a curve.

Despite its usefulness, there are few criteria that help in
computing the rank of maps among various spaces of sections of
vector bundles (see \cite{B}, \cite{T2}). The purpose of the
following note is to help fill in this gap. We want to show that
for a {\bf generic} curve and {\bf any} line bundle of degree at
most $g+1$ with a given number of sections, the symmetric product
map
$$S^2(H^0(C,L))\rightarrow H^0(C,L^2)$$
is injective. The result is obviously false for line bundles of
high degree, as even the dimension of the left hand side space
becomes bigger than the dimension of the right hand side. We prove
in fact that our bound is sharp: we construct an example of a line
bundle of degree $g+2$ on a generic curve of  genus $g$ (for every
even $g\ge 4$) for which the symmetric product map is not
injective (see \ref{exemple}).

Note that the natural product map
$$m:\  H^0(C,L)\otimes H^0(C,L)\rightarrow H^0(C,L^2)$$
is never injective for $k\ge 2$. Its kernel contains the subspace
$\wedge ^2(H^0(C,L))$. Therefore, proving the injectivity of the
map with domain in the symmetric power is equivalent to showing
that the kernel of $m$ is precisely the wedge power. This kernel
has proved  to have  important deformation theoretical meaning in
particular cases (see for instance [T1]).

The main result of this paper is the following:

\begin{Thm} \label{teorema} {\bf Theorem}
Let $C$ be a generic curve defined over an algebraically closed
field of characteristic different from two. Let $L$ be a line
bundle on $C$ such that $deg(L)=d\le g+1$. Then the map
$$S^2(H^0(C,L))\rightarrow H^0(C,L^2)$$ is injective.
\end{Thm}

The proof of this result is inspired in the proof of the
injectivity of the Petri map given by Eisenbud and Harris in
\cite{EH} (see also \cite{W}).

\end{section}

\begin{section}{The case of degree at most $g-1$}

 We shall fix a genus $g$ and the degree $d$ for line bundles $L$ on $C$.
 We shall denote by $k$ (rather than the classical
$r+1$) the dimension of the space of sections of these line
bundles. We shall assume that $k, g, d$ have been fixed so that
the generic curve of genus $g$ has line bundles of degree $d$ with
$k$ sections (equivalently, the Brill- Noether number
$g-k(g-1-d+k)$ is non-negative).

 In order to prove \ref{teorema}  for a generic curve, it suffices to
prove it for a special curve. Consider a family of curves $\pi
:{\cal C}\rightarrow T$. Let $T$ be the spectrum of a discrete
valuation ring ${\cal O}$ with maximal ideal generated by $t$.
Assume that the generic fiber of $\pi$ is a non-singular curve and
the special fiber $C$ looks as follows:

Take $g$ elliptic curves $E^i$ and let $P^i,Q^i$ be generic points
on $E^i$. Take any number of rational curves
$C^0_1,..C^0_{k_0},...C^g_0...C^g_{k_g}$ again with points $P^i_j,
Q^i_j$ on them. Glue $C^i_j$ to $C^i_{j+1}$ by identifying $Q^i_j$
to $P^i_{j+1}$. Glue $C^{i-1}_{k_{i-1}}$ to $E^i$ by identifying
$Q^{i-1}_{k_{i-1}}$ to $P^i$. Glue $E^i$ to $C^i_1$ by identifying
$Q^i$ to $P^i_1$.

 For convenience of notation, we shall
denote by
$$Y_1,...Y_M,\ M=k_0+...+k_g+g$$  the components of $C$ starting with $C^0_1$
and ending with $C^g_{k_g}$. We shall denote by $P_i,Q_i$ the two
points in $Y_i$ that get identified to $Q_{i-1}\in Y_{i-1}$ and
$P_{i+1}\in Y_{i+1}$ respectively. We warn the reader that we
shall keep the superindices $i$ when we need to refer to the
$i^{th}$ elliptic curve. We hope this will not produce too much
confusion.

Note that the form of the central fiber does not change if we make
base changes and normalisations.

Consider now the set up of  limit of a linear series  in the sense
of Eisenbud and Harris (\cite{EH} p.273 or \cite{EH1}, section 2).
Let ${\cal L}$ be a line bundle on ${\cal C}$. One can modify
${\cal L}$ by tensoring with divisor with support on the central
fiber. This leaves invariant the line bundle on the generic fiber
but modifies it in the central fiber. For every component $Y_i$ of
$C$, there is a line bundle ${\cal L}_i$ on ${\cal C}$ such that
it has degree zero on every component of the central fiber except
for the component $Y_i$. As ${\cal L}_i= {\cal L}_{i+1}(-d\sum
_{j\le i}Y_j)$, one can identify ${\cal L}_i$ with a subsheaf of
${\cal L}_{i+1}$.

 Consider
$\pi_*({\cal L}_i)$ This is a free ${\cal O}$ module of rank {k}.
Moreover, ${\cal L}_i\subset {\cal L}_{i+1}$ is a lattice. Denote
by $V_i$ the image of the restriction map
$$\pi _*{\cal L}_i\rightarrow \pi _*{\cal L}_{i|Y_i}=H^0(Y_i,{\cal L}_{i|Y_i})$$
As $deg{\cal L}_{i|Y_j}=0,\ j\not= i$, this map is injective and
we shall sometimes identify $\pi_*({\cal L}_i)$ with $V_i$.

From \cite{EH} Lemma 1.2, one can find a basis of sections
$\sigma_m^i, m=1..k$ of the free module $\pi_*({\cal L}_i)$ such
that the orders of vanishing of the sections $\sigma_m^i$ at $P_i$
are the different orders of vanishing of the sections of $V_i$ and
$t^{\alpha _m^i}\sigma _m^i,\ m=1...k$ form a basis for $\pi
_*{\cal L}_{i+1}$.

We shall now relate the  vanishing of sections of line bundles at
the various nodes.

\begin{Thm} {\bf Lemma} \label{Lemma1} 1) Let $Y$ be any irreducible non-singular curve  $L$ a
line bundle of degree $d$ on $Y$ and $P,Q$ two points on $Y$. The
sum of the orders of vanishing at $P$ and $Q$ of any section of
$L$ is at most $d$.

2)  Let $Y$ be an elliptic curve and $P,Q$ generic points of $Y$.
Let $L$ be a line bundle of degree $d$ on $Y$.The sum of the
orders of vanishing at $P$ and $Q$ of any section of $L$ is at
most $d-1$ except in the case where $L={\cal O}(aP+(d-a)Q)$ for
some $a$. In this case, there is only one section of $L$ vanishing
with multiplicities adding up to $d$ at the two points.
\end{Thm}

\begin{pf} If a section of a line bundle $L$ vanishes with order $a$
at $P$ and $b$ at $Q$, then ${\cal O}(aP+bQ)$ is a subsheaf of
$L$. Hence, $d\ge a+b$. This proves the first statement.

If $a+b=d$, then $L={\cal O}(aP+bQ)$. If there is another section
vanishing to orders $a',b'$ at $P,Q$ with $a'+b'=d$ and say $a'>a$
then $aP+bQ\equiv a'P+b'Q$. Hence, $cP\equiv cQ$ with $c=a'-a$.
This contradicts the genericity of the pair $P,Q$ if $Y$ is not
rational.
\end{pf}

{\bf Remark} The genericity of $P,Q$ is essential here. If $cP$ is
linearly equivalent to $cQ$ for some $c\le d$, then the line
bundle ${\cal O}(aP+(d-a)Q), a\ge c$ has (at least) two sections
with orders of vanishing adding up to $d$ at $P,Q$ namely
$aP+(d-a)Q$ and $(a-c)P+(d-a+c)Q$.
\bigskip

The following result of Eisenbud and Harris (cf. Prop 1.1 in
[E,H]), will be used in the sequel.

\begin{Thm}\label{anul} {\bf Lemma} Let $\sigma$ be a section in $\pi _*{\cal L}_i$.
Let $\alpha$ be the unique integer such that $t^{\alpha}\sigma \in
\pi_*({\cal L}_{i+1})-t\pi _*({\cal L}_{i+1})$, then
$$ord_{P_i}(\sigma _{|Y_i})\le d-ord_{Q_i}(\sigma _{|Y_i})\le \alpha \le
ord_{P_{i+1}}(t^{\alpha}\sigma _{|Y_{i+1}})$$
\end{Thm}

As in [E,H], p.277, one can define the order of  a section $\rho
\in S^2(\pi_*{\cal L}_Y)$ at a point $P$ in a component $Y$ as
follows:

\begin{Thm} {\bf Definition} We say $ord_P(\rho _{|Y})\ge l$ if and only if $\rho$ is in the
linear span of $t(S^2(\pi_*({\cal L}_Y)))$ and elements of the
form $\sigma_m\otimes \sigma_n+\sigma_n\otimes \sigma_m$ where
$ord_P(\sigma_n)+ord_P(\sigma _m)\ge l,\ \sigma_n,\ \sigma_m \in
\pi_*({\cal L}_Y)$.
\end{Thm}

One then has the following result (cf. \cite{EH}, Lemma 3.2)

\begin{Thm}\label{anulrho} {\bf Lemma} Let $\sigma_m$ be a basis
of the free ${\cal O}$ module $\pi_*({\cal L}_i)$ such that the
orders of vanishing of the $\sigma _m$ at $P_i$ are the distinct
orders of vanishing of the linear series at this point and
$t^{\alpha _m}\sigma _m$ is a basis of $\pi_*({\cal L}_{i+1})$. If
$$\rho =\sum f_{n,m}(\sigma_n\otimes\sigma_m+\sigma_m\otimes
\sigma_n)\in S^2(\pi_*{\cal L}_i)$$ where the $f_{n,m}$ are
functions on the discrete valuation ring ${\cal O}$ and the
associated discrete valuation is $\nu$, then
$$ord_{P_i}(\rho_{|Y_i})=min_{\{\nu(f_{n,m})=0\} } (ord_{P_i}(\sigma_n)+ord_{P_i}(\sigma_m))$$
If $\beta$ is the unique integer such that
$$t^{\beta}\rho \in S^2\pi_*{\cal L}_{i+1}-t(S^2\pi_*{\cal
L}_{i+1})$$ then
$$\beta =max\{ \alpha _n+\alpha _m-\nu (f_{nm})\}$$
\end{Thm}

Let us assume now that the kernel of the symmetric product map is
non-zero on the generic curve. We can then find an element $\rho$
such that say $\rho \in S^2(\pi_*{\cal L}_1)-tS^2(\pi_*{\cal
L}_1)$ and  integers $\beta_i, i=2...M$ such that $t^{\beta
_i}\rho \in S^2(\pi_*{\cal L}_i)-tS^2(\pi_*{\cal L}_i)$ and
$t^{\beta _i}\rho \in Ker(S^2(\pi_*{\cal L}_i)\rightarrow
\pi_*({\cal L}_i)^2)$.

As a section of a line bundle cannot vanish to order higher than
the degree, the following claim will conclude the proof in the
case $d\le g-1$.

\begin{Thm} {\bf Claim} \label{claim} If $l> k_0+...+k_{m-1}+m$, then $ord_{P_l}(t^{\beta _l} \rho
)\ge 2m$. In particular, for $l\ge k_0+...+k_{g-1} +g,\
ord_{P_l}(t^{\beta _l} \rho )\ge 2g$.
\end{Thm}

\begin{pf} We prove the
following two statements:

1) If $C_i$ is a rational curve,
$$ord_{P_{i+1}}(t^{\beta_{i+1}}\rho_{|Y_{i+1}})\ge
ord_{P_i}(t^{\beta_i}\rho_{|{Y_i}})$$.

2) If $C_i$ is an elliptic component,
$$ord_{P_{i+1}}(t^{\beta_{i+1}}\rho_{|Y_{i+1}})\ge
ord_{P_i}(t^{\beta_i}\rho_{|{Y_i}})+2$$

In other words, the order of vanishing of a section
$t^{\beta_{i+1}}\rho$ at $P_{i+1}$ is at least as large as that of
$t^{\beta_i}\rho$ at $P_i$ if $Y_i$ is rational and at least two
units larger if it is elliptic. As the order of vanishing
$t^{\beta_0}\rho$ at $P_0$ is non-negative, this shows that the
order of vanishing of $t^{\beta_i}\rho$ at $P_i$ is at least $2m$
if $m$ elliptic curves precede $Y_i$. This is the first part of
the claim. The second part follows from the first when $m=g$.

We now turn to the proof of 1) and 2).  Choose a basis $\sigma
_m,\ m=1...k$ of $\pi_*({\cal L}_i)$ such that $t^{\alpha
_m}\sigma _m$ is a basis of $\pi_*{\cal L} _{i+1}$. For simplicity
of notation, we shall assume that $\beta_i=0$ Write
$$\rho =\sum _{n\le m}f_{nm}(\sigma_n\otimes \sigma _m +\sigma
_m\otimes \sigma_n)\in S^2(\pi_* {\cal L}_i)$$ Then, from
\ref{anulrho},
$$ord_{P_i}(\rho _{|Y_i})=min_{\{ \nu
(f_{nm})=0\}}(ord_{P_i}(\sigma_n)+ord_{P_i}(\sigma_m))$$ Assume
that this minimum is attained by a pair corresponding to the
indices  $n_0,m_0$ with $\nu(f_{n_0,m_0})=0$. Then from
\ref{anul},
$$(ord_{P_i}(\sigma_{n_0})+ord_{P_i}(\sigma_{m_0}))\le
2d- ord_{Q_i}(\sigma_{n_0})-ord_{Q_i}(\sigma_{m_0}) \le
\alpha_{n_0}+\alpha_{m_0}$$ From \ref{anulrho} and the fact that
$\nu (f_{n_0,m_0})=0$, the latter is at most $\beta_{i+1}$.

Write $$t^{\beta_{i+1}}\rho =\sum _{n\le
m}(t^{\beta_{i+1}-\alpha_n-\alpha_m}f_{nm})(t^{\alpha_n}\sigma_n\otimes
t^{\alpha_m}\sigma _m +t^{\alpha_m}\sigma _m\otimes
t^{\alpha_n}\sigma_n)$$

 Hence, from \ref{anulrho} $$ord_{P_{i+1}}(t^{\beta_{i+1}}\rho
_{|Y_{i+1}})=min_{\{ \beta_{i+1}-\alpha_n-\alpha_m+\nu
(f_{nm})=0\}}(ord_{P_{i+1}}(t^{\alpha_n}\sigma_n)+ord_{P_{i+1}}(t^{\alpha
_m}\sigma_m))$$

Assume that this minimum is attained at a pair $n_1,m_1$ with
$$\beta_{i+1}-\alpha_{n_1}-\alpha_{m_1}+\nu(f_{n_1,m_1})=0$$
Then, $$\beta_{i+1}\le
\beta_{i+1}+\nu(f_{n_1,m_1})=\alpha_{n_1}+\alpha_{m_1}\le
ord_{P_{i+1}}(t^{\alpha_{n_1}}\sigma_{n_1})+ord_{P_{i+1}}(t^{\alpha
_{m_1}}\sigma_{m_1}))$$ where the last inequality comes from
\ref{anul}

 Stringing together the above inequalities, we obtain
$ord_{P_i}(\rho_{|Y_i})\le ord_{P_{i+1}}(\rho_{|Y_{i+1}})$. Hence
part 1) is proved.

 Assume now that there is equality in the
inequality above. Then all the previous inequalities must be
equalities. In particular, the terms $\sigma_{n_0}\otimes
\sigma_{m_0}+\sigma_{m_0}\otimes \sigma_{n_0}$ that give the
vanishing of $\rho$ at $P_i$ satisfy
$$ord_{P_i}(\sigma _k)+ord_{Q_i}(\sigma_k)=d,\ k=n_0,m_0$$.

If we have
$ord_{P_{i+1}}(\rho_{|Y_{i+1}})=ord_{P_i}(\rho_{|Y_i})+1$, then in
each pair, one of $\sigma_k,\ k=n_0, m_0$ would vanish to order
$d$ between the two nodes. From \ref{Lemma1}, there is at most one
such section $\sigma _{i_0}$ on our elliptic curves.

Hence, if $ord_{P_{i+1}}(\rho_{|Y_{i+1}})\le
ord_{P_i}(\rho_{|Y_i})+1$, the terms in $\rho$ giving the
vanishing at $P_i$ could be written as
$$\sigma_{i_0}\otimes \sigma +\sigma\otimes \sigma_{i_0}$$
for some section $\sigma$. As a section like this cannot be in the
kernel of the symmetric product map, this concludes the  proof of
the claim.

\end{pf}
\end{section}
\begin{section}{The cases of degree $g, g+1$ and counterexamples
in degree $g+2$}

Assume now $d=g$. Consider the last elliptic component $E^g$ with
nodes $P^g,Q^g$. Write the restriction to the central fiber of an
element in the kernel of the symmetric evaluation map as
$$\rho=\sum f_{nm}(\sigma_n\otimes \sigma_m+\sigma _m\otimes \sigma_n),
\ \sigma_n,\ \sigma_m\in \pi_*({\cal L}_{E^g} )$$.

Consider a pair $\sigma_n,\sigma_m$ that gives the vanishing of
$\rho$ at $P^g$. From \ref{claim},
$ord_{P^g}(\sigma_n)+ord_{P^g}(\sigma_m)\ge 2(g-1)$. As $deg({\cal
L}_{E^g|E^g})=g$, $ord_{P^g}(\sigma_k)\le g,\ k=n,m$. It follows
that the possible orders of such sections are $g, g-1, g-2$. If
there is a section vanishing to order $g$ at $P^g$, then ${\cal
L}_{E^g|E^g}={\cal O}(gP^g)$ and this line bundle does not have
any section vanishing to order $g-1$ at $P^g$. Hence there are at
most two independent sections among the $\sigma_k$ and for every
pair that appears, the sum of the orders of vanishing at $P^g$ is
the same. But then $\rho$ cannot be in the kernel of the symmetric
product map. This concludes the proof for $d=g$

Assume now that $d=g+1$. If $ord_{P^1}(\rho)=0,1$, then the pairs
of sections giving the vanishing at $P^1$ will vanish to orders
$0,1$. From our choice of basis, there are at most two sections
satisfying these conditions. Hence, again $\rho$ could not be in
the kernel of the symmetric product map. Therefore the vanishing
of $\rho$ at $P^1$ is at least two. From statements 1) and 2) in
\ref{claim}, the vanishing at $P^g$ is at least $2g$. Now, we
conclude as in the case $d=g$.
\bigskip

\begin{Thm} \label{exemple} {\bf Example} If $g$ is even,
$g\ge 4$ a generic curve of genus $g$ has a line bundle of degree
$g+2$ such that the kernel of the symmetric product map is not
injective.
\end{Thm}
\begin{pf}  A generic curve of even genus $g$ has line bundles
of degree $g/2+1$ with two sections (because the corresponding
Brill-Noether number is zero). Take two different ones $L_1, L_2$
and choose independent sections $s_1,t_1$ and $s_2, t_2$ of them.
Take $L=L_1\otimes L_2$. This is a line bundle of degree $g+2$
that contains the sections $s_1s_2, s_1t_2,t_1s_2, t_1t_2$. These
are four independent sections, as can be checked using the base
point free pencil trick.

 Consider now the following section
$$(s_1s_2\otimes t_1t_2+t_1t_2\otimes s_1s_2)-(s_1t_2\otimes
s_2t_1+s_2t_1\otimes s_1t_2)$$ This section is non-zero and in the
kernel of the symmetric product map.
\end{pf}

\end{section}

\end{document}